\documentclass[11pt]{amsart}

\usepackage{amsfonts}
\usepackage{amsthm}
\usepackage{amssymb}
\usepackage{amsmath}
\usepackage{url, hyperref}

\newtheorem{theorem}{Theorem}
\newtheorem{lemma}{Lemma}
\newtheorem{corollary}{Corollary}
\newtheorem{remark}{Remark}
\newtheorem{proposition}{Proposition}

\newcommand{\Z}{{\mathbb Z}}

\newcommand{\R}{{\mathbb R}}

\newcommand{\bea}{\begin{eqnarray*}}
\newcommand{\eea}{\end{eqnarray*}}
\newcommand{\be}{\begin{eqnarray}}
\newcommand{\ee}{\end{eqnarray}}

\newcommand{\prob}{\mbox{\rm Prob}\,}
\newcommand{\frob}{\mbox{\rm F}\,}

\numberwithin{equation}{section}

\begin{document}

\title[Expected Frobenius numbers]{Expected Frobenius numbers}
\author{Iskander Aliev}
\address{School of Mathematics and Wales Institute of Mathematical and Computational Sciences, Cardiff University, Senghennydd Road, Cardiff, Wales, UK}
\email{alievi@cf.ac.uk}

\author{Martin Henk}
\address{Institut f\"ur Algebra und Geometrie, Otto-von-Guericke
Universit\"at Mag\-deburg, Universit\"atsplatz 2, D-39106-Magdeburg}
\email{martin.henk@ovgu.de}

\author{Aicke Hinrichs}
\address{Friedrich Schiller Universit\"at Jena\\
Fakult{\"a}t für Mathematik und Informatik\\
Ernst-Abbe-Platz 2\\
D-07743 Jena}
\email{hinrichs@minet.uni-jena.de}


\begin{abstract}
We show that for large instances the order of magnitude of the expected Frobenius number is (up to a constant depending only on the dimension) given by its lower bound.
\end{abstract}

\maketitle

\section{Introduction}
Let $a$ be a positive integral $n$-dimensional primitive vector,  i.e., $a=(a_1,\dots,a_n)^\intercal\in\Z_{>0}^n$ with $\gcd(a):=\gcd(a_1,\dots,a_n)=1$. The {\em Frobenius number} of $a$, denoted by $\frob(a)$, is the largest number which cannot be represented as a non-negative integral combination of the $a_i$'s, i.e.,
\begin{equation*}
  \frob(a)=\max\{b\in\Z : b\neq \langle a, z\rangle \text{ for all } z\in \Z_{\geq 0}^n\},
\end{equation*}
where $\langle \cdot, \cdot\rangle$ denotes the standard inner product on $\R^n$.
 In other words, $\frob(a)$ is the maximal right hand side $b$, such that the well-known {\em knapsack polytope} $P(a,b)=\{x\in\R^n_{\geq 0} : \langle a, x\rangle=b\}$ does not contain an integral point. From that point of view it is also apparent that the Frobenius number plays an important role in the analysis of integer programming algorithms (see, e.g., \cite{Aardal:2002p6506, Hansen_Ryan, Lee:2008p7096, Scarf1993, Sturmfels:1996p7884}) and, vice versa, integer programming algorithms are known to be an effective tool for computing the Frobenius number  (see, e.g., \cite{Beihoffer:2005p6531, Einstein:2007p7808, Roune:2008p7660}).
There is a rich literature on Frobenius numbers, and for an impressive survey on the history and the different aspects of the problem  we refer to the book \cite{RamirezAlfonsin:2005p6858}.

Here we just want to mention that only for $n=2$ an explicit formula is known, which was proven by  Curran Sharp in 1884 answering a question by Sylvester:
\begin{equation*}
       \frob(a)=a_1\,a_2-(a_1+a_2).
\end{equation*}
There is a huge variety of upper bounds on $\frob(a)$. They all share the property
that in the worst case they are of quadratic order with respect to the  maximum norm of $a$, say, which will be denoted by $|a|_\infty$. For instance, assuming $a_1\le a_2\le \ldots\le a_n$, a classical upper
bound due to  Erd\H{o}s and Graham \cite{ErdHos:1972p6594} says
\begin{equation*}
\frob(a)\le 2\, a_n\,\left [\frac{a_1}{n}\right ]-a_1,
\end{equation*}
and, in a recent paper, Fukshansky\&Robins \cite[(29)]{Fukshansky:2007p6629} gave an upper bound which is also symmetric in the $a_i$'s
\begin{equation}
 \frob(a)\leq \left[\frac{(n-1)^2/\Gamma(\frac{n}{2}+1)}{\pi^{n/2}}\sum_{i=1}^na_i\sqrt{(\vert a\vert_2)^2-a_i^2}+1\right],
\label{eq:fukshansky_robins}
\end{equation}
where $|\cdot|_2$ denotes the Euclidean norm. The worst case in all the known upper bounds is achieved when the $a_i$'s are approximately of the same size, and it is also known that in these cases the quadratic order of an upper bound cannot be lowered (see \cite{Arnold:2006p151, ErdHos:1972p6594, SchlagePuchta:2005p7045}).

On the other hand, Aliev\&Gruber \cite{Aliev:2007p6944} recently found an optimal lower bound for the Frobenius number which implies that 
\begin{equation}
\frob(a)> (n-1)!^{\frac{1}{n-1}}\,\left(a_1\,a_2\cdot\ldots\cdot a_n\right)^{\frac{1}{n-1}}- (a_1+\cdots +a_n)\,.
\label{eq:aliev_gruber}
\end{equation}
 Hence, if all the $a_i$'s are of the same size then the lower bound is only of order $|a|_\infty^{1+1/(n-1)}$. In fact, taking the quotient of  the (symmetric) upper bound \eqref{eq:fukshansky_robins} with \eqref{eq:aliev_gruber}, we see that  there is always a gap of order $|a|_\infty^{1-1/(n-1)}$.

Thus the next natural and important question is to get information on the Frobenius number of a ``typical'' vector $a$. This problem appears to be hard, and to the best of our knowledge it has firstly been systematically investigated by V.~I.~Arnold, see, e.g., \cite{Arnold:1999p6499, Arnold:2006p151, Arnold:2007p7106}.
In particular, he conjectured that $\frob(a)$ grows like $T^{1+1/(n-1)}$ for a ``typical'' vector $a$ with  $1$-norm $\vert a\vert_1=T$ \cite{Arnold:1999p6499}, and in  \cite[2003-5]{Arnold:2004p6389} he conjectures that the ''average behavior'' is
\begin{equation}
     \frob(a)\sim (n-1)!^{\frac{1}{n-1}}\,\left(a_1\,a_2\cdot\ldots\cdot a_n\right)^{\frac{1}{n-1}},
\label{eq:conj_arnold}
\end{equation}
i.e., it is essentially the lower bound. A similar conjecture for the $3$-dimen\-sional case was proposed  by Davison \cite{Davison:1994p6498}  and recently proved by Shur, Sinai and Ustinov \cite{Shur:2008p2269}.  Extensive computations  support conjecture \eqref{eq:conj_arnold} (see \cite{Beihoffer:2005p6531}).

In \cite{Bourgain:2007p560},  Bourgain and  Sinai proved a statement in the spirit of these conjectures, which says, roughly speaking, that
\begin{equation*}
 \prob\left(\frob(a)/T^{1+1/(n-1)}\geq D \right) \leq \epsilon(D),
\label{eq:BS}
\end{equation*}
where $\prob(\cdot)$ is meant with respect to the uniform
distribution among all points in the set
\begin{equation*}
   G(T)=\{ a\in\Z^n_{>0}: \gcd(a)=1,\,\vert a\vert_\infty\leq T\}.
\end{equation*}
Here the number $\epsilon(D)$ does not depend on $T$ and tends to zero as $D$ approaches infinity. The paper \cite{AlievHenk} gives more precise information about the order of decay of the function $\epsilon(D)$. Their main result \cite[Theorem 1.1]{AlievHenk} implies that
\begin{equation}
 \prob\left(\frob(a)/\vert a\vert_\infty^{1+1/(n-1)}\geq D \right) \ll_n D^{-2}\,,
\label{eq:alievhenk_frob}
\end{equation}
where $\ll_n$ denotes the Vinogradov symbol with the constant depending on $n$ only.
 In particular, from that result the authors get a statement about the average Frobenius number, namely \cite[Corollary 1.2]{AlievHenk}
\begin{equation*}
   \sup_T  \frac{\sum_{a\in \mathrm{G}(T)}\frob(a)/\vert a\vert_\infty^{1+1/(n-1)}}{\#\mathrm{G}(T) }\ll_n 1.
\end{equation*}
So the order of the average  Frobenius number is (essentially) not bigger than $\vert a\vert_\infty^{1+1/(n-1)}$, which is close to the sharp lower bound \eqref{eq:aliev_gruber}, but there is still a gap.

The main purpose of this note is to fill that gap. We will show that
\begin{theorem} Let $n\geq 3$. Then
\begin{equation*}
 \prob\left(\frob(a)/\left(a_1\,a_2\cdot\ldots\cdot a_n\right)^{\frac{1}{n-1}} \geq D \right) \ll_n D^{-2\frac{n-1}{n+1}}.
\end{equation*}
\label{thm:main}
\end{theorem}
From this result  we will derive the desired statement
\begin{corollary} Let $n\geq 3$. Then
\begin{equation*}
 \sup_T  \frac{\sum_{a\in \mathrm{G}(T)}\frob(a)/\left(a_1\,a_2\cdot\ldots\cdot a_n\right)^{\frac{1}{n-1}}}{\#\mathrm{G}(T) }\ll \gg_n1.
\end{equation*}
\label{cor:main}
\end{corollary}
These statements fit also perfectly to recent results on the limit distribution of Frobenius numbers due to Shur, Sinai and Ustinov \cite{Shur:2008p2269} and Marklof \cite{Marklof:2009p2303}. For instance, in our special setting, \cite[Theorem 1]{Marklof:2009p2303} says that
\begin{equation*}
\lim_{T\to\infty} \prob\left(\frob(a)/\left(a_1\,a_2\cdot\ldots\cdot a_n\right)^{\frac{1}{n-1}} \geq  D \right) =\Psi(D),
\end{equation*}
where $\Psi:\R_{\geq 0}\to \R_{\geq 0}$ is a non-increasing function with $\Psi(0)=1$.

The proof of Theorem \ref{thm:main} is based on a discrete inverse arithmetic-geometric mean inequality which might be of some interests in its own. It will be stated and proved in Section 2. Finally, Section 3 contains the proofs of Theorem \ref{thm:main} and Corollary \ref{cor:main}.

\section{Reverse discrete AGM-inequality}
For $x\in\R^n_{\geq 0}$ the Arithmetic-Geometric-Mean (AGM) inequality states that $\left(\prod_{i=1}^n x_i\right)^{1/n}\leq\frac{1}{n}\sum_{i=1}^n x_i$. It is known that the ``reverse'' AGM inequality holds with high probability. More precisely, Gluskin\&Milman \cite{Gluskin:2003p3223} have shown that
\begin{equation*}
 \prob_{S^{n-1}}\left(\frac{\sqrt{\frac{1}{n}\sum_{i=1}^n  x_i^2}}{\left(\prod_{i=1}^n \vert x_i\vert\right)^{1/n}}=\frac{1}{\sqrt{n}\left(\prod_{i=1}^n \vert x_i\vert\right)^{1/n}} >\alpha\right )\leq c^n\,\alpha^{-n/2},
\end{equation*}
where $c$ is an absolute constant and $ \prob_{S^{n-1}}()$ is meant with respect to standard rotation invariant measure on the $n$-dimensional unit sphere $S^{n-1}=\{x\in\R^n : \sum_{i=1}^n x_i^2 =1\}$. Here we show an analogous statement with respect to the primitive lattice points in the set $G(T)$.

With respect to the uniform distribution on $G(T)$, let $L_T:G(T)\to\R_{> 0}$ be the random variable defined by
\begin{equation*}
         L_T(a) = \frac{\frac{1}{n}\sum_{i=1}^n  a_i}{\left(\prod_{i=1}^n  a_i\right)^{1/n}}.
\end{equation*}
\begin{theorem} Let $\alpha > 1$ and let $k\in\{1,\dots,n-1\}$. Then there exists a constant $c(k,n)$ depending only on $k$, $n$, such that
\begin{equation*}
   \prob(L_T\geq \alpha)\leq c(k,n)\,\alpha^{-k}.
\end{equation*}
\label{thm:reverse_agm}
\end{theorem}
By Markov's inequality the theorem is an immediate consequence of the next statement about the expectation $\mathbb{E}( L_T ^k)$ of higher moments of $L_T$.
\begin{lemma} Let $k\in\{1,\dots,n-1\}$. Then there exists a constant
$c(k,n)$ depending on $k$, $n$, such that
\begin{equation*}
   \mathbb{E}( L_T ^k)\leq c(k,n).
\end{equation*}
\label{lem:expectation}
\end{lemma}
\begin{proof}  First we note that there is an absolute  constant $c$ such that
\begin{equation}
          \# G(T)\geq c\,T^n.
\label{eq:counting_primitive}
\end{equation}
This follows easily from well-known relations between integration and counting primitive lattice points (see e.g. \cite[p. 183, (1)]{Cassels:1971p7266}), but
in order to keep the paper self-contained we give a short argument here: Let $n\geq 2$ and let 
$$G_2(T)=\{ (a,b)^\intercal \in \Z_{\geq 1}^2 : 1\le a,b \le T, \gcd(a,b)=1\}.$$
Since $G_2(T)\times \{1,\dots,T\}^{n-2}\subseteq G(T)$ it suffices to prove the statement for $n=2$, i.e., $G_2(T)$. There are at most $(T/m)^2$ pairs $(a,b)^\intercal \in \{1,\dots,T\}$  with $\gcd(a,b)=m$. Thus
\begin{equation*}
  \# G_2(T)\geq T^2\left(1-\sum_{m=2}^\infty m^{-2}\right)= (2-\pi^2/6)\,T^2,
\end{equation*}
which gives \eqref{eq:counting_primitive}.

Now we have
\begin{equation}
\begin{split}
 \mathbb{E}( L_T ^k) & = \frac{1}{ \# G(T)} \sum_{a\in G(T)} L_T(a)^k = \frac{1}{ \# G(T)}\sum_{a\in G(T)}\left( \frac{\frac{1}{n}\sum_{i=1}^n a_i}{\left(\prod_{i=1}^n a_i\right)^\frac{1}{n}}\right)^k \\
&=\frac{1}{ \# G(T)}\sum_{a\in G(T)} \frac{1}{n^k}\sum_{i_1+i_2+\cdots + i_n=k}\binom{k}{i_1,i_2,\cdots,i_n} \frac{a_1^{i_1}\,a_2^{i_2}\cdot\ldots\cdot a_n^{i_n}}{\left(\prod_{i=1}^n a_i\right)^\frac{k}{n}}\\
&\leq \frac{1}{ \# G(T)}\sum_{i_1+i_2+\cdots + i_n=k}\frac{1}{n^k}\binom{k}{i_1,i_2,\cdots,i_n} \sum_{a_1,a_2,\ldots,a_n=1}^T \prod_{j=1}^n a_j^{i_j-k/n} \\
&=\frac{1}{ \# G(T)}\sum_{i_1+i_2+\cdots + i_n=k}\frac{1}{n^k}\binom{k}{i_1,i_2,\cdots,i_n}  \prod_{j=1}^n \sum_{a_j=1}^T a_j^{i_j-k/n}.
\end{split}
\label{eq:long}
\end{equation}
Since for $k<n$ the sum $\sum_{m=1}^T m^{-k/n}$ is bounded from above by $c(k,n)\,T^{1-k/n}$, where $c(k,n)$ is a constant depending only on $k$ and $n$, we find
\begin{equation}
   \sum_{a_j=1}^T a_j^{i_j-k/n} \leq c(k,n) T^{1+i_j-k/n}.
\label{eq:sum_int}
\end{equation}
Thus, for $i_1+i_2+\cdots + i_n=k$ we obtain
\begin{equation*}
       \prod_{j=1}^n \sum_{a_j=1}^T a_j^{i_j-k/n} \leq c(k,n)^n\,T^n.
\end{equation*}
Hence we can continue  \eqref{eq:long} by
\begin{equation*}
\begin{split}
   \mathbb{E}( L_T ^k) & \leq \frac{1}{ \# G(T)}\sum_{i_1+i_2+\cdots + i_n=k}\frac{1}{n^k}\binom{k}{i_1,i_2,\cdots,i_n}\, c(k,n)^n\,T^n\\ &= c(k,n)^n\,\frac{T^n}{ \# G(T)}.
\end{split}
\end{equation*}
Finally, with \eqref{eq:counting_primitive} we get the assertion.
\end{proof}

\begin{remark} With a little more work one can
prove Lemma \ref{lem:expectation} for any positive number $k<n$, not only for integers.
\label{rem:case3}
\end{remark}

Next we want to point out that the arguments in the proof of Lemma \ref{lem:expectation} lead to the following lower bound on the random variable
$X_T:G(T)\to \R_{>0}$ given by
\begin{equation*}
         X_T(a) = \frac{\frob(a)}{\left(\prod_{i=1}^n  a_i\right)^{\frac{1}{n-1}}}.
\end{equation*}

\begin{proposition} Let $n\geq 2$. Then there exists a constant depending only on $n$ such that
\begin{equation*}
        \mathbb{E}(X_T)\geq (n-1)!^\frac{1}{n-1}\left(1-c(n)\,T^{-\frac{1}{n-1}}\right).
\end{equation*}
\label{prop:lower_expectation}
\end{proposition}
\begin{proof} On account of the lower bound \eqref{eq:aliev_gruber} on $\frob(a)$ it remains to show that
\begin{equation*}
             \frac{1}{\# G(T)}\sum_{a\in G(T)} \frac{\sum_{i=1}^na_i}{\left(\prod_{i=1}^n a_i\right)^\frac{1}{n-1}}\leq c(n)\,T^{-\frac{1}{n-1}}.
\end{equation*}
Following the argumentation in \eqref{eq:long} we find
\begin{equation*}
\begin{split}
\sum_{a\in G(T)} \frac{\sum_{i=1}^na_i}{\left(\prod_{i=1}^n a_i\right)^\frac{1}{n-1}} &\leq
\sum_{i_1+i_2+\cdots + i_n=1} \sum_{a_1,a_2,\ldots,a_n=1}^T \prod_{j=1}^n a_j^{i_j-\frac{1}{n-1}}\\
&= n\, \sum_{a_1,a_2,\ldots,a_n=1}^T a_1^{{1-\frac{1}{n-1}}}a_2^{-\frac{1}{n-1}}\cdot\ldots\cdot a_n^{-\frac{1}{n-1}}\\
&=n\,\sum_{a_1=1}^T a_1^{{1-\frac{1}{n-1}}}\left(\sum_{a_2=1}^T a_2^{{-\frac{1}{n-1}}}\right)^{n-1}.
\end{split}
\end{equation*}
Thus, analogously to \eqref{eq:sum_int}, and on account of \eqref{eq:counting_primitive} we obtain
\begin{equation*}
\begin{split}
\frac{1}{\# G(T)}\sum_{a\in G(T)} \frac{\sum_{i=1}^na_i}{\left(\prod_{i=1}^n a_i\right)^\frac{1}{n-1}} &\leq \tilde{c}(n)\,
\frac{1}{\# G(T)} T^{2-\frac{1}{n-1}}\left(T ^{1-\frac{1}{n-1}}\right)^{n-1}\\ &\leq c(n) T^{-\frac{1}{n-1}}.
\end{split}
\end{equation*}
\end{proof}

\section{Proof of Theorem \ref{thm:main} and Corollary \ref{cor:main}}
We keep the notation of the previous section.
First we note  that \eqref{eq:alievhenk_frob} is certainly also true for any other norm 
in the denominator, in particular for the 1-norm $|\cdot|_1$. Thus
\begin{equation}
 \prob\left(\frac{\frob(a)}{\vert a\vert_1^{1+\frac{1}{n-1}}}\geq D \right) \ll_n D^{-2}.
\label{eq:step1}
\end{equation}
Secondly, we observe that
\begin{equation}
\begin{split}
\prob\left(\frac{|a|_1^{1+\frac{1}{n-1}}}{(a_1\cdot\ldots\cdot a_n)^{\frac{1}{n-1}}}\geq \gamma \right) &=\prob\left(\frac{(\frac{1}{n}|a|_1)^{\frac{n}{n-1}}}{(a_1\cdot\ldots\cdot a_n)^{\frac{1}{n-1}}}\geq \left(\frac{1}{n}\right)^{\frac{n}{n-1}}\gamma \right)\\
&=\prob\left(L_T^\frac{n}{n-1}\geq \left(\frac{1}{n}\right)^{\frac{n}{n-1}}\gamma \right)\\
&=\prob\left(L_T\geq \frac{1}{n}\gamma^\frac{n-1}{n} \right)\ll_n \gamma^{-\frac{(n-1)^2}{n}},
\end{split}
\end{equation}
where in the last inequality we have applied Theorem \ref{thm:reverse_agm} with $k=n-1$. Together with \eqref{eq:step1} we obtain
\begin{equation*}
\begin{split}
\prob\left( X_T(a)\geq \beta\right) & = \prob\left(\frac{\frob(a)}{\vert a\vert_1^{1+\frac{1}{n-1}}}\cdot \frac{|a|_1^{1+\frac{1}{n-1}}}{(a_1\cdot\ldots\cdot a_n)^{\frac{1}{n-1}}}\geq \beta\right)\\
&\leq \prob\left(\frac{\frob(a)}{\vert a\vert_1^{1+\frac{1}{n-1}}} \geq \beta^t\right) +
\prob\left(\frac{|a|_1^{1+\frac{1}{n-1}}}{(a_1\cdot\ldots\cdot a_n)^{\frac{1}{n-1}}}\geq \beta^{1-t}\right)\\
&\ll_n \beta^{-2\,t} + \beta^{-\frac{(n-1)^2}{n} (1-t)},
\end{split}
\end{equation*}
for any $t\in(0,1)$. With $t=(n-1)/(n+1)$ we finally get
\begin{equation*}
\prob\left( X_T(a)\geq \beta\right) \ll_n \beta^{-2\frac{n-1}{n+1}},
\end{equation*}
which shows Theorem \ref{thm:main}.

For the proof of Corollary \ref{cor:main} we note that 
\begin{equation*}
\begin{split}
\mathbb{E}(X_T)&=\int_{0}^\infty \prob(X_T>x)\,\mathrm{d}x\leq 1+\int_{1}^\infty \prob(X_T>x)\,\mathrm{d}x\\
&\ll_n 1+\int_{1}^\infty x^{-2\frac{n-1}{n+1}}\,\mathrm{d}x.
\end{split}
\end{equation*}
For $n\geq 4$, the last integral is finite and so we have
\begin{equation}
\mathbb{E}(X_T) \ll_n 1.
\label{eq:upper_expectation}
\end{equation}
For the case $n=3$ we just note that on account of  Remark \ref{rem:case3}
one can also bound $\prob(X_T(a)\geq \beta)$ by a function like $\beta^{-(1+\epsilon)}$
and so we also get \eqref{eq:upper_expectation} in this case. Together with Proposition \ref{prop:lower_expectation}, Corollary \ref{cor:main} is proven.


\begin{thebibliography}{99}

\bibitem{Aardal:2002p6506}
K. Aardal and A. K. Lenstra, \emph{Hard equality constrained integer
  knapsacks}, Lect. Notes Comput. Sci. \textbf{2337} (2002), 350--366.

\bibitem{RamirezAlfonsin:2005p6858}
J.~L~Ram{\'\i}rez Alfons{\'\i}n, \emph{The Diophantine Frobenius problem},
  Oxford Lecture Series in Mathematics and its Applications \textbf{30} (2005),
  xvi+243.

\bibitem{Aliev:2007p6944}
I.~M. Aliev and P.~M. Gruber, \emph{An optimal lower bound for the
  Frobenius problem}, Journal of Number Theory \textbf{123} (2007), no.~1,
  71--79.

\bibitem{AlievHenk}
I.~M. Aliev and M. Henk, \emph{Integer knapsacks: Average behavior of the Frobenius numbers},  to appear in Mathematics of Operations Research, arxiv:0810.0234v1.

\bibitem{Arnold:1999p6499}
V. I. Arnold, \emph{Weak asymptotics for the numbers of solutions of Diophantine
  problems}, Functional Analysis and Its Applications \textbf{33} (1999),
  no.~4, 292--293.

\bibitem{Arnold:2004p6389}
\bysame, \emph{Arnold's problems}, Springer-Verlag, Berlin (2004).

\bibitem{Arnold:2006p151}
\bysame, \emph{Geometry and growth rate of frobenius numbers of additive
  semigroups}, Math Phys Anal Geom \textbf{9} (2006), no.~2, 95--108.

\bibitem{Arnold:2007p7106}
\bysame, \emph{Arithmetical turbulence of selfsimilar fluctuations statistics
  of large Frobenius numbers of additive semigroups of integers}, Mosc. Math.
  J. \textbf{7} (2007), no.~2, 173--193, 349.




\bibitem{Beihoffer:2005p6531}
D. Beihoffer, J. Hendry, A. Nijenhuis, and S. Wagon, \emph{Faster
  algorithms for Frobenius numbers}, Electron. J. Combin. \textbf{12} (2005),
  Research Paper 27, 38 pp. (electronic).

\bibitem{Bourgain:2007p560}
J. Bourgain and Y.~G. Sinai, \emph{Limit behaviour of large Frobenius
  numbers}, Russ. Math. Surv. \textbf{62} (2007), no.~4, 713--725.

\bibitem{Cassels:1971p7266}
J. W. S Cassels, \emph{An introduction to the Geometry of Numbers},
 {Grundlehren der mathematischen Wissenschaften}, (1971), {Springer-Verlag, Berlin-New York}.

\bibitem{Davison:1994p6498}
J.L Davison, \emph{On the linear Diophantine problem of Frobenius}, Journal of
  Number Theory \textbf{48} (1994), 353--364.

\bibitem{Einstein:2007p7808}
D. Einstein, D. Lichtblau, A. Strzebonski and St. Wagon,
\emph{Frobenius numbers by lattice point enumeration},
Integers \textbf{7} (2007), 63 pages.


\bibitem{ErdHos:1972p6594}
P~Erd\H{o}s and R.~L Graham, \emph{On a linear Diophantine problem of Frobenius},
  Acta Arith. \textbf{21} (1972), 399--408.

\bibitem{Fukshansky:2007p6629}
L. Fukshansky and S. Robins, \emph{Frobenius problem and the covering
  radius of a lattice}, Discrete Comput. Geom. \textbf{37} (2007), no.~3,
  471--483.



\bibitem{Gluskin:2003p3223}
E.~Gluskin and V.~Milman, \emph{Note on the geometric-arithmetic mean inequalty},
  Lecture Notes in Mathematics \textbf{1807} (2003), 131--135.

\bibitem{Hansen_Ryan}
P. Hansen, J. Ryan,
\emph{Testing integer knapsacks for feasibility},
European Journal of Operational Research \textbf{88} (1996), no. 3, 578--582.


\bibitem{Lee:2008p7096}
J. Lee, S. Onn, and R. Weismantel, \emph{On test sets for nonlinear
  integer maximization}, Oper. Res. Lett. \textbf{36} (2008), no.~4, 439--443.

\bibitem{Marklof:2009p2303}
J. Marklof, \emph{The asymptotic distribution of frobenius numbers},  arXiv:0902.3557,  (2009), 14 pages.

\bibitem{Roune:2008p7660}
B. Hammersholt Roune, 
\emph{Solving thousand-digit Frobenius problems using Gr{\"o}bner bases}, 
J. Symbolic Comput. \textbf{43} (2008), no. 1, 1--7.

\bibitem{Scarf1993}
H.~E. Scarf and D.~F. Shallcross, \emph{The frobenius problem and maximal
  lattice free bodies}, Math. Oper. Res. \textbf{18} (1993), no.~3, 511 -- 515.

\bibitem{SchlagePuchta:2005p7045}
J.-C. Schlage-Puchta, \emph{An estimate for frobenius' diophantine
  problem in three dimensions}, J. Integer Seq. \textbf{8} (2005), no.~1,
  Article 05.1.7, 4 pp. (electronic).

\bibitem{Shur:2008p2269}
V.~Shur, Y.~Sinai, and A.~Ustinov, \emph{Limiting distribution of Frobenius
  numbers for $n=3$}, arXiv:0810.5219,  (2008), 13 pages.


\bibitem{Sturmfels:1996p7884}
B. Sturmfels, 
\emph{Gr{\"o}bner bases and convex polytopes},
University Lecture Series 8 (1996), 
 American Mathematical Society, Providence.
\end{thebibliography}

\providecommand{\bysame}{\leavevmode\hbox to3em{\hrulefill}\thinspace}
\providecommand{\MR}{\relax\ifhmode\unskip\space\fi MR }
\providecommand{\MRhref}[2]{%
  \href{http://www.ams.org/mathscinet-getitem?mr=#1}{#2}
}
\providecommand{\href}[2]{#2}

\end{document}